\documentclass{article}

{}
{}
{}
\usepackage{color}
\usepackage{multirow}
\usepackage{makecell}
\usepackage{float}
\usepackage{setspace}
\usepackage{array}
\usepackage{longtable}
\usepackage{supertabular}
\usepackage{geometry}
\usepackage{fancyhdr} 

\usepackage{mathtools}
\usepackage{amssymb}
\usepackage{graphicx}
\usepackage{lineno,hyperref}

\usepackage{amsfonts}
\usepackage{comment}

\begin{document}

\title{Estimating Lower Probability Bound of Power System's Capability to Fully Accommodate Variable Wind Generation}

\author{Bin Liu\thanks{School of Electrical Engineering and Telecommunications, The University of New South Wales, Sydney 2052, Australia. Email address: eeliubin@hotmail.com.}, 
	Bingxu Zhai\thanks{Jibei Electric Power Dispatching and Control Center, State Grid Cooperation of China, Beijing 100054, China.},
	Mengchen Liu\thanks{Department of Grid Connection, Goldwind Australia Pty Ltd, Sydney 2000, Australia.},
	Feng Liu\thanks{Department of Electrical Engineering, Tsinghua University, Beijing 100084, China.},
 	Haibo Lan\thanks{Jibei Electric Power Dispatching and Control Center, State Grid Cooperation of China, Beijing 100054, China.}}

\maketitle
\thispagestyle{fancy}          
\fancyhead{}                      
\lhead{This paper is a preprint of a paper accepted by The Journal of Engineering and is subject to Institution of Engineering and Technology Copyright. When the final version is published, the copy of record will be available at the IET Digital Library.}           
\chead{}
\rhead{}
\lfoot{}
\cfoot{}   
\rfoot{\thepage}
\renewcommand{\headrulewidth}{0pt}       
\renewcommand{\footrulewidth}{0pt}

\begin{abstract}
As the penetration of wind generation increases, the uncertainty it brings has imposed great challenges to power system operation. To cope with the challenges, tremendous research work has been conducted, among which two aspects are of most importance, i.e. making immune operation strategies and accessing the power system's capability to accommodate the variable energy. Driven and inspired by the latter problem, this paper will discuss the power system's capability to accommodate variable wind generation in a probability sense. Wind generation, along with its uncertainty is illustrated by a polyhedron, which contains prediction, risk and uncertainty information. Then, a three-level optimization problem is presented to estimate the lower probability bound of power system's capability to fully accommodate wind generation. After reformulating the inner \emph{max-min} problem, or feasibility check problem, into its equivalent mixed-integer linear program (MILP) form, the bisection algorithm is presented to solve this challenging problem. Modified IEEE systems are adopted to show the effectiveness of the proposed method. 
\end{abstract}

\section{Introduction}
Developing renewable energies, especially wind generation, has been an efficient and effective way for many countries in recent years to solve energy  deficiency problem and  environmental issues. However, wind generation are weak in controllability, which is well equipped in traditional thermal or hydro generators \cite{ref1,ref27}. The resulted wind generation uncertainty has imposed great challenge on power system operation like unit commitment (UC) and economic dispatch (ED) \cite{ref13,ref14,ref15}. To cope with this problem, tremendous work has been conducted either by exploring the flexibility of current power system or integrating other devices with quicker response like energy storage systems (EES) \cite{ees-1,ees-2,ees-3}. In terms of improving power system flexibility, improving the accuracy of modeling wind generation uncertainty, making immune operation strategies and evaluating power system's capability to accommodate wind generation are three topics that have been wildly concerned. 

Generally, wind generation uncertainty can be depicted as a a number of scenarios associated with certain occurrence probabilities \cite{ref2,ref3} or uncertainty sets \cite{ref4,ref5,ref6} which can be modeled as a ellipsoid or polyhedron. Then, stochastic UC or ED is developed based on wind scenarios and the optimal solution is feasible for all scenarios while minimizing the expected power operation cost \cite{ref7,ref8,ref9,ref24}. Similarly, robust UC or ED is formulated based on uncertainty set of wind generation and the optimal solution is feasible for any possible scenario belonging to the uncertainty set \cite{ref10,ref11,ref12}. According to different mathematical forms, the robust UC or ED either minimize the power operation cost under predicted wind generation or under worst wind scenario belonging to the constructed uncertainty set. On the contrary, given a dispatching strategy, the power system can be assessed on how much wind generation it can accommodate, which leads to the problem to be discussed in this paper. For example, once the unit commitment is determined either by deterministic, stochastic or robust method, a set of wind scenarios can be generated to test the effectiveness of the solution by simulating a series of ED problem. The percentage of wind scenarios associated with a feasible ED problem can be used to depict the effectiveness of the given strategy. Actually, the widely used Monte Carlo method to evaluate the effectiveness of novel dispatching strategy in research articles can be regarded as a method of assessing the capability of power system to accommodate wind generation or coping with wind generation uncertainty based on generated wind scenarios \cite{ref16}. The concept of dispatchable region of variable wind generation is presented as a polyhedron by solving a series of linear program (LP) and mixed-integer linear program (MILP) problems in \cite{dregion}. The dispatchable region is also a measure of flexibility of given strategies. In \cite{windrisk}, a risk-based admissibility assessment approach is proposed to evaluate how much wind generation can be accommodated after presenting a measure on the risk that variable wind generation bring to power system.  

Driven by the idea of detecting infeasibility in robust optimization, we present a novel method to estimate the lower probability bound of power system's capability to \emph{{fully}} accommodate variable wind generation. Different from existing works, in this paper, after bridging a given probability with the presented uncertainty set in polyhedron form, the estimation model is formulated to maximize this probability, which provides a metric to depict the guaranteed chance that power system can fully accommodate variable wind generation. The estimation model is formulated as a three-level or \emph{max-max-min} optimization problem. Then, a bisection algorithm is presented after reformulating the inner \emph{max-min} problem to its equivalent MILP form. 

The paper is organized as follows. The polyhedral uncertainty set of wind generation, which is associated with a probability, is presented in Part 2. The UC model, which provides unit commitments for estimation model is given in Part 3. In Part 4, the estimation model is formulated, followed by the bisection algorithm. Case study is performed in Part 5 and the paper is concluded in Part 6.

\section{Constructing Uncertainty Set}
In robust optimization, uncertainty set is generally formulated as a closed set to depict the possible realizations of wind generation. The uncertainty set is usually constructed beforehand based on physical nature or historical data of wind generation and then put into the robust optimization as know quantity. In this part, the uncertainty set will be constructed based on the statistic characteristics of historical data, as we explain next.

Denote the possible wind generation as $w$, with its expectation value and prediction error as $w_e$ and $\Delta w=w-w_e$, respectively. The covariance matrix of $\Delta w$, denoted by $\Lambda$, can be obtained from analysis of historical data. Then, we have the following equation.
\begin{equation}\label{eq-1}
w-w_e=\Lambda^{1/2}v~~\mathbf{or}~~v=\Lambda^{-1/2}(w-w_e)
\end{equation}
where $w=[w_{1,1},\cdots,w_{i,h},\cdots,w_{W,H}]^T$ with $W,H$ representing the total number of wind farm and dispatch interval, respectively; Denoting $m=W\times H$, then $v\in \mathbb{R}^{m\times 1}$ is a vectorial variable, whose expected value and covariance matrix are $\mathbf{0}_{m\times 1}$ and $\mathbf{1}_{m\times m}$, respectively; $\mathbb{R}$ is the set of real numbers.

Obviously, the elements in $v$ are independent to each other, making it much easier and more straightforward to construct the uncertainty set by imposing constraints on $v$. Besides, bridging the probability and the uncertainty set will also be showed as follows.   

Assuming the probability of $v_i$ falling into the box $[-z_i,z_i]$ is $\gamma_i=\varphi(z_i)$, i.e. $\varphi(z_i)=\textbf{Prob}\{-z_i\le v_i\le z_i\}$, and let $\prod_i\gamma_i=\alpha$. Construct the following uncertainty set.
\begin{eqnarray} 
\label{eq-2-1}
\mathcal{W}=\{w|-u\le v=\Lambda^{-1/2}(w-w_e)\le u\} \\
\label{eq-2-2}
u=[\varphi^{-1}(\gamma_1),\cdots,\varphi^{-1}(\gamma_m)]^T
\end{eqnarray} 

Then, the possible wind generation falling into the uncertainty set $\mathcal{W}$ is $\alpha$. \begin{eqnarray} 
\label{eq-2-3}
\textbf{Prob}\{w\in\mathcal{W}\}=\alpha
\end{eqnarray} 

Then, the only question to be addressed is how to determine $\gamma_i~(\forall i)$ with given $\alpha$. Considering that all elements in $v$ are independent to each other, we simply assume that $\gamma_i~(\forall i)$ are set to the same value in this paper, resulting to $\gamma_i=\sqrt[m]{\alpha}~(\forall i)$. Therefore, given the probability of $\alpha$, the uncertainty set $\mathcal{W}$ can be constructed as
\begin{equation}\label{eq-3}
\begin{split}
\mathcal{W(\alpha)}=\{w|-\varphi^{-1}(\sqrt[m]{\alpha})&\le\Lambda^{-1/2}(w-w_e)\le\varphi^{-1}(\sqrt[m]{\alpha})\} \\
\end{split}
\end{equation}

Besides, we simply assume wind prediction error is normally distributed in this paper, which has been widely studied in literatures. This means $v$ is also normally distributed according to \eqref{eq-1}. However, other kinds of distribution pattern can also be employed \cite{error-1,error-2,error-3}.

\section{Unit Commitment Formulation}
The UC problem is used to provide the unit commitments of thermal generators, which are part of known information in the estimation model. As the UC problem can also be helpful to understanding the formulation of the estimation model, it is presented as follows. 
\begin{eqnarray} 
\label{uc-1}
F_{UC}=\min\sum_{i=1}^{N}\sum_{h=1}^{H}{S_iv_{i,h}}+\sum_{i=1}^{N}\sum_{h=1}^{H}{(B_ip_{i,h}+C_iu_{i,h})}\\
\label{uc-2}
s.t.~~v_{i,h}\ge u_{i,h}-u_{i,h-1}~\forall i,\forall h\\
\label{uc-3}
u_{i,k}\ge u_{i,h}-u_{i,h-1}\\
~\forall i,\forall h, \forall k\in\{h,\cdots,h+T_i^{\text{on}}-1\}\nonumber\\
\label{uc-4}
u_{i,k}\le 1-u_{i,h-1}+u_{i,h}\\
~\forall i,\forall h, \forall k\in\{h,\cdots,h+T_i^{\text{off}}-1\}\nonumber\\
\label{uc-5}
u_{i,h}P_i^-\le p_{i,h}\le u_{i,h}P_i^+~\forall i,\forall h\\
\label{uc-6}
p_{i,h+1}-p_{i,h}\le u_{i,h}R_i^++(1-u_{i,h})P_i^+~\forall i,\forall h\\
\label{uc-7}
p_{i,h}-p_{i,h+1}\le u_{i,h+1}R_i^-+(1-u_{i,h+1})P_i^+~\forall i,\forall h\\
\label{uc-8}
\sum_{i=1}^{N}{p_{i,h}}+\sum_{i=1}^{N}{w_{i,h}}=\sum_{i=1}^{N}{D_{i,h}}~\forall h\\
\label{uc-9}
\sum_{i=1}^{N}{(u_{i,h}P_i^+-p_{i,h})}\ge r_h^+\sum_{i=1}^{N}{D_{i,h}}~\forall h\\
\label{uc-10}
\sum_{i=1}^{N}{(p_{i,h}-u_{i,h}P_i^-)}\ge r_h^-\sum_{i=1}^{N}{D_{i,h}}~\forall h\\
\label{uc-11}
-L_k\le\sum_{i=1}^{N}{H_{i,k}(p_{i,h}+w_{i,h}-D_{i,h})}\le L_k~\forall k, \forall h
\end{eqnarray}
where, $u_{i,h}$ is a binary variable indicating the on/off state of thermal generator $i$ at hour $h$, i.e. $u_{i,h}=1$ when the generator $i$ is 'ON' at hour $h$ and $u_{i,h}=0$ otherwise; $v_{i,h}$ is also a binary variable, which indicates whether the generator $i$ is started up from hour $h$, i.e. $v_{i,h}=1$ means generator $i$ starts up in hour $h$ and $v_{i,h}=0$ otherwise; $p_{i,h}$ is continuous variable representing the generation of thermal generator $i$ at hour $h$; $S_i,B_i,C_i$ are start-up cost, generation cost per unit and minimum operation cost of thermal generator $i$, respectively; $T_i^{\text{on}}/T_i^{\text{off}}$ is the minimum on/off hours of thermal generator $i$; $P_i^-/P_i^+$ is the lower/upper generation bound of thermal generator $i$;  $R_i^{-}/R_i^{+}$ is ramping down/up rate of thermal generator $i$; $D_{i,h}$ is the active power load of bus $i$ at hour $h$; $r_h^{-}/r_h^{+}$ is down/up reserve requirement of power system at hour $h$; $L_k$ is the active power limit of transmission line $k$ and $H_{i,k}$ is the sensitivity of power injection of bus $i$ to power flow in transmission line $k$.

In UC model, the objective is formulated as \eqref{uc-1} and \eqref{uc-2}-\eqref{uc-4} makes sure that the start-up is captured and the minimum on/off hour requirement is met for any thermal generator. Equation \eqref{uc-5} guarantees that generation is within the limit of each thermal generator and \eqref{uc-6}-\eqref{uc-7} means the thermal generation can only be adjusted within the ramping capability. Equation \eqref{uc-8} makes sure the active power is balanced in each hour while \eqref{uc-9}-\eqref{uc-10} guarantees that the up/down reserve is sufficient in each hour. Finally, the limits of transmitting active power by each line is imposed by \eqref{uc-11}.

The optimal solution $u^*=\{u_{i,h}~\forall i, \forall h\}$ of UC will be provided to the estimation model presented in the next part.

\section{Estimation Method Formulation}
\subsection{\it Problem Formulation}
As discussed previously, the obtained optimal solution $u^*$ will be input as known parameter in the estimation model. In other words, the estimation model will evaluate how much wind generation can be fully accommodated under the UC strategy  $u^*$. Before formulating the estimation model, we present the following abbreviations.
\begin{eqnarray} 
p_{i,h}^-=u^*_{i,h}P_i^-,~p_{i,h}^+=u^*_{i,h}P_i^+\nonumber\\
\Delta p_{i,h}^-=u^*_{i,h+1}R_i^-+(1-u^*_{i,h+1})P_i^+\nonumber\\
\Delta p_{i,h}^+=u^*_{i,h}R_i^++(1-u^*_{i,h})P_i^+\nonumber\\
\Delta r_h^-=\sum_{i=1}^{N}{u^*_{i,h}P_i^-}+r_h^-\sum_{i=1}^{N}{D_{i,h}}\nonumber\\
\Delta r_h^+=\sum_{i=1}^{N}{u^*_{i,h}P_i^+}-r_h^+\sum_{i=1}^{N}{D_{i,h}}\nonumber
\end{eqnarray}

For each $w\in\mathcal{W(\alpha)}$, the feasibility check problem can be formulated as ($\textbf{CK}$ for short)
\begin{eqnarray} 
\label{ck-1}
F_{c}(w)=\min\sum_{i=1}^{N}\sum_{h=1}^{H}{(z_{i,h}^{a+}}+z_{i,h}^{a-})+\sum_{h=1}^{H-1}{(z_{i,h}^{b+}}+z_{i,h}^{b-})\\
+\sum_{h=1}^{H}{(z_{h}^{c+}+z_{h}^{c-}+z_{h}^{d+}+z_{h}^{d-})}+\sum_{k=1}^{K}\sum_{h=1}^{H}{(z_{k,h}^{e+}+z_{k,h}^{e-})}\nonumber\\
\label{ck-2}
s.t.~~~~~~~~~~~~~~~~~~~~~~~~p_i^--z_{i,h}^{a-}\le p_{i,h}~\forall i,\forall h\\
\label{ck-3}
p_{i,h}-z_{i,h}^{a+}\le p_i^+~\forall i,\forall h\\
\label{ck-4}
p_{i,h+1}-p_{i,h}-z_{i,h}^{b+}\le \Delta p_{i,h}^+~\forall i,\forall h\\
\label{ck-5}
p_{i,h}-p_{i,h+1}-z_{i,h}^{b-}\le \Delta p_{i,h}^-~\forall i,\forall h\\
\label{ck-6}
\sum_{i=1}^{N}{p_{i,h}}+\sum_{i=1}^{N}{w_{i,h}}+z_{h}^{c+}-z_{h}^{c-}=\sum_{i=1}^{N}{D_{i,h}}~\forall h\\
\label{ck-7}
\Delta r_h^--z_{h}^{d-}\le \sum_{i=1}^{N}{p_{i,h}}~\forall h\\
\label{ck-8}
\sum_{i=1}^{N}{p_{i,h}}-z_{h}^{d+}\le \Delta r_h^+~\forall h\\
\label{ck-9}
-L_k-z_{k,h}^{e-}\le\sum_{i=1}^{N}{H_{i,k}(p_{i,h}+w_{i,h}-D_{i,h})}~\forall k, \forall h\\
\label{ck-10}
\sum_{i=1}^{N}{H_{i,k}(p_{i,h}+w_{i,h}-D_{i,h})}-z_{k,h}^{e+}\le L_k~\forall k, \forall h\\
\label{ck-11}
z_{i,h}^{a+},z_{i,h}^{a-},z_{i,h}^{b+},z_{i,h}^{b-},z_{k,h}^{e+},z_{k,h}^{e-}\ge 0,~\forall i, \forall k, \forall h\\
\label{ck-12}
z_{h}^{c+},z_{h}^{c-},z_{h}^{d+},z_{h}^{d-}\ge 0,~\forall h
\end{eqnarray}
where $z_{i,h}^{a+},z_{i,h}^{a-},z_{i,h}^{b+},z_{i,h}^{b-},z_{h}^{c+},z_{h}^{c-},z_{h}^{d+},z_{h}^{d-},z_{k,h}^{e+},z_{k,h}^{e-}$ are introduced slack variables.  

In $\textbf{CK}$, a positive optimal objective value, i.e. $F_c(w)>0$, means the optimal value of at least one slack variable is positive, indicating that no feasible dispatching strategy for $w$ is available in real-time dispatch under $u^*$.

Based on $\textbf{CK}$, the estimation problem can be formulated as ($\textbf{AP}$ for short) 
\begin{eqnarray} 
\label{asse-1}
&\alpha_{0}=\max{\alpha}\\
&s.t.~~~~~F_{w}(\alpha)= 0\\
\label{asse-2}
&F_{w}(\alpha) = \max\limits_{w\in\mathcal{W}(\alpha)}\min\limits_{p,z}\left\{ 1^Tz\Bigg|
\begin{array}{ll}
Mp+Tw+Qz\le m,~(\delta)\\
z\ge 0\\
\end{array} \right\}
\end{eqnarray}
where the inner \emph{min} problem of \eqref{asse-2} is the compact form of $\textbf{CK}$; $M,Q,T$ are known matrices of appropriate dimension and $p,w,z$ represent the vector of thermal generation, wind generation and slack variable;  $\delta$ is the Lagrange multiplier of the corresponding inequality. 

Similar to the explanation of $F_c(w)$, $F_{w}(\alpha)>0$ means at least one possible wind generation in $\mathcal{W}(\alpha)$ is unable to be balanced by thermal generators in real-time operation under $u^*$. On the contrary, $F_{w}(\alpha)=0$ means any realization of $w\in\mathcal{W}(\alpha)$ can be balanced in real-time operation. As $\alpha$ is used to quantify the possibility of wind generation falling into $\mathcal{W}(\alpha)$, the maximum value of $\alpha$ that satisfies $F_{w}(\alpha)=0$, i.e. $\alpha_0$, actually provides a lower bound of power system's probability, or guaranteed chance, to {\it fully} accommodate variable wind generation. In other words, $1-\alpha$ indicates the upper bound of the chance that  curtailing wind generation occurs in the real-time operation under $u^*$. 

In practical operation , the information obtained can be provided to power system dispatcher to help them understand the system's capability to accommodate wind generation better, as well as help them adjust UC strategy or make necessary preparations.

\subsection{\it Algorithm}

Obviously, $\textbf{AP}$ is a complicated three-level or \emph{max-max-min} problem. To solve this  challenging problem, we first reformulate the inner \emph{max-min} problem to its equivalent MILP form, as we explain next.

By dualizing the innermost \emph{min} problem, the \emph{max-min} problem \eqref{asse-2} can be reformulated as ($\textbf{CK}_{dual}$ for short) 
\begin{eqnarray} 
\label{dual-1}
F_{w}(\alpha) =\max\limits_{w\in\mathcal{W}(\alpha)}\max\limits_{\delta}{\delta^TTw-\delta^Tm}\\
\label{dual-2}
s.t.~~~~~~1+Q^T\delta\ge 0\\
\label{dual-3}
M^T\delta=0,~\delta\ge 0
\end{eqnarray}

Noting that the optimal value of $\textbf{CK}_{dual}$ must be achieved at the extreme point of $\mathcal{W}(\alpha)$, it can be further reformulated as the following MILP problem ($\textbf{CK}_{milp}$ for short).
\begin{eqnarray} 
\label{milp-1}
F_{w}(\alpha) =\max\limits_{t,\tau,\rho,\delta}{\sum\limits_{n}{(t_n+\rho_nx_n^-)}-\delta^Tm}\\
\label{milp-2}
s.t.~~~~~~1+Q^T\delta\ge 0\\
\label{milp-3}
M^T\delta=0,~\delta\ge 0\\
\label{milp-4}
t_n\le\tau_n\rho_n^+(x_n^+-x_n^-)~\forall n\\
\label{milp-5}
t_n\le[\rho_n-(1-\tau_n)\rho_n^-](x_n^+-x_n^-)~\forall n\\
\label{milp-6}
\tau_n\in\{0,1\}~\forall n
\end{eqnarray}
where $\rho=\Lambda^{1/2}T^T\delta$, $x=\Lambda^{-1/2}w$ with $x^+=\Lambda^{-1/2}w_e+u$, $x^-=\Lambda^{-1/2}w_e-u$ and $\rho^+/\rho^-$ is the estimated upper/lower bound of $\rho$.

$\textbf{CK}_{milp}$ can be efficiently solved by commercial solvers such as Cplex \cite{cplex} and Gurobi \cite{gurobi}. By solving a series of MILP problem, $\alpha_0$ can be obtained by the following bisection algorithm.  

\textbf{Algorithm 1}

\textbf{Step 1}: Initialize error tolerance $\varepsilon_{\alpha},\varepsilon_{F}$, upper and lower bounds of $\alpha$ as $\alpha^+,\alpha^-$.

\textbf{Step 2}: If $|\alpha^+-\alpha^-|\le\varepsilon_{\alpha}$, go to \textbf{Step 4}. Otherwise, go to the next step.

\textbf{Step 3}: Set $\alpha=\frac{\alpha^++\alpha^-}{2}$ and solve $\textbf{CK}_{milp}$. Denote the optimal value as $F_w^*$.

1) If $F_w^*\le\varepsilon_{F}$, update $\alpha^-=\frac{\alpha^++\alpha^-}{2}$, go to \textbf{Step 2}.

2) If $F_w^*>\varepsilon_{F}$, update $\alpha^+=\frac{\alpha^++\alpha^-}{2}$, go to \textbf{Step 2}.

\textbf{Step 4}: Record the optimal solution of $\textbf{AP}$, i.e. $\alpha_0=\frac{\alpha^++\alpha^-}{2}$.

In this paper, all programs are coded with YALMIP \cite{yalmip} in Matlab and all MILP problems are solved by CPLEX in a Thinkpad with Intel(R) i7-3520M 2.9 GHz and 8 GB RAM.
\section{Case Study}
\subsection{Setup of the tested system}
\begin{figure}[H]
	\centering\includegraphics[scale=0.25]{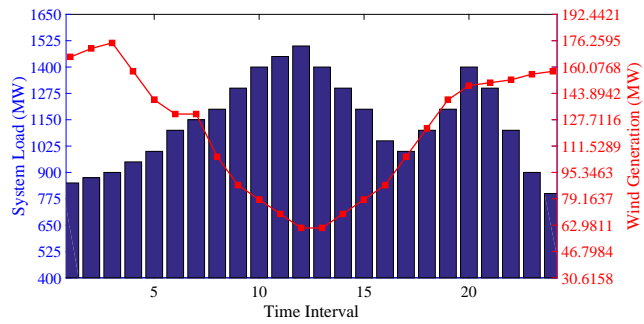}
	\caption{Predicted wind generation \& system load}
	\label{fig-1}
\end{figure}

In this section, the modified IEEE 39-bus system is studied to demonstrate the effectiveness of the proposed method. In the modified case, one wind farm is connected to bus 10 and the prediction of wind generation, along with the system load is given in Fig.\ref{fig-1}. 
The wind generation profile is obtained from the historical data of EIRGRID \cite{ref26}, however, with adjustment to fit the system load level. Besides, we assume that the prediction information of wind generation in the following day is updated before 0:00, which leads to the assumption that the accuracy of wind generation prediction declines form 0:00 to 24:00 of the next day. Accordingly, we set 
\begin{eqnarray} 
\label{wind-1}
\sigma_{i}=0.05+0.05\times (i-1)/23\\
\rho_{i,j}= \left\{ \begin{array}{ll}
1 & \textrm{if $i=j$}\\
\rho_{i,i}=1-0.9\times |j-i|/23 & \textrm{if $i\neq j$}\\
\end{array} \right.
\end{eqnarray} 
where $\sigma_{i}$ is the standard error of $\Delta \hat{w}_i=\Delta w_i/w_e$, which is the relative prediction error of wind generation; $\rho_{i,j}$ is the correlation coefficient between  $\Delta \hat{w}_i$ and  $\Delta \hat{w}_j$.

The minimum up and down reserve requirements both are 5\% of power system load in each hour, i.e. $r^-_h=r^+_h=5\%~(\forall h)$. Other parameters of thermal generators are showed in Table \ref{tab-1}. The positive/negative integer in the 'Ini.' and 'Loc.' columns indicate the hours of on/off state thermal generator has stayed at the beginning hour and the locations of thermal generators. 
%
\begin{table}[H]
	\setlength{\tabcolsep}{1.5pt}
	\renewcommand\arraystretch{0.9}
	\caption{Parameters of thermal generators}
	\begin{tabular}{>{\small}c>{\small}c>{\small}c>{\small}c>{\small}c>{\small}c>{\small}c>{\small}c>{\small}c}
		\hline
		No.&Loc.& Ini.& \makecell{$P_i^-/P_i^+$\\(MW)}&\makecell{$R_i^-/R_i^+$\\(MW/h)}&\makecell{$S_i$\\(\$)} &\makecell{$T^{\text{on}}/T^{\text{off}}$\\(h)} &	\makecell{$B_i$\\(\$/MWh)}&	\makecell{$C_i$\\(\$)}\\
		\hline
		1	&30&	8	&	150/455 &	70/70	&	6750	&	8/8	&	16.19	&	1000\\
		2	&31&	8	&	150/455 &	60/60	&	7500	&	8/8	&	17.26	&	970	\\
		3	&32&	-5	&	20/130	&	70/70	&	825	    &	5/5	&	16.60	&	700	\\
		4	&33&	-5	&	20/130 	&	70/70	&	840  	&	5/5	&	16.50	&	680	\\
		5	&34&	-6	&	25/162	&	60/60	&	1350	&	6/6	&	19.70	&	450	\\
		6	&35&	-3	&	20/80	&	70/70	&	255	    &	3/3	&	22.26	&	370	\\
		7	&36&	-3	&	25/85	&	70/70	&	390	    &	3/3	&	27.74	&	480	\\
		8	&37&	-1	&	10/55	&	60/60	&	45	    &	1/1	&	25.92	&	660	\\
		9	&38&	-1	&	10/55	&	60/60	&	45	    &	1/1	&	27.27	&	665	\\
		10	&39&	-1	&	10/55	&	60/60	&	45	    &	1/1	&	27.79	&	670	\\
		\hline
	\end{tabular}\centering
	\label{tab-1}\end{table}
\subsection{Uncertainty set construction}
\begin{figure}[H]
	\centering\includegraphics[scale=0.25]{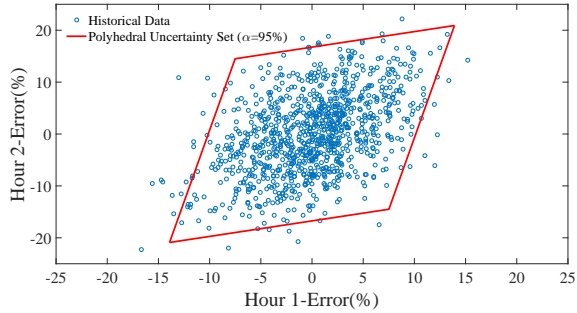}
	\caption{Constructing polyhedral uncertainty set}
	\label{fig-2}
\end{figure}

To show the effectiveness of uncertainty set constructed in Part 2, the simulation data of hour 1 and hour 15 for the connected wind farm is studied. According to previous introduction, the covariance matrix of $\Delta \hat{w}_1$ and $\Delta \hat{w}_{15}$ is [0.00250, 0.00182; 0.00182, 0.00647]. The uncertainty set is constructed in Fig.\ref{fig-2}. It is noteworthy that the presented uncertainty set in Fig.\ref{fig-2} is for relative prediction error and the uncertainty set for wind generation can be constructed by moving the center to $w_e$.
In the above figure, the uncertainty set is constructed with $\alpha=95.00\%$. Then 2000 prediction error scenarios are randomly generated purely according to the distribution information wind prediction error. The simulation was repeated for 100 times and the average percentage of points falling into the uncertainty set is 95.05\%, which is very close to $\alpha$. Such a simulation result demonstrates that the polyhedral uncertainty set depicts the uncertainty characteristics of prediction error effectively. 
\subsection{Estimation results}
\begin{figure}[H]
	\centering\includegraphics[scale=0.25]{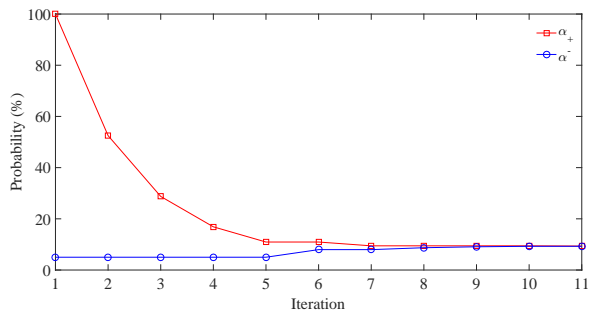}
	\caption{Iteration information of solving $\textbf{AP}$}
	\label{fig-3}
\end{figure}

The UC problem is solved firstly and the optimal operation strategy $u^*$ is then put into the estimation model as known parameters. After 36.73 seconds and 11 iterations, the optimal solution is obtained as $\alpha_0$=9.36\%, which means the guaranteed chance that the power system can fully accommodate variable generation under $u^*$ is 9.36\%. The iteration process, demonstrated by $\alpha^+$ and $\alpha^-$, are also showed in in Fig. \ref{fig-3}. 
The result also indicates that there is a chance of up to 90.64\% that curtailing wind generation will occur in real-time operation of the coming day. As $u^*$ is made by the deterministic UC model presented in Part 3, the result also implies that making a UC decision purely through minimizing the operation cost may leads to a strategy that is not the best one to accommodate wind generation, and that employing novel methods, e.g. stochastic and robust UC, to make the decision is necessary when facing the uncertainty wind generation brings to power system.
\section{Conclusion}
In this paper, we discussed the estimation of the lower probability bound of power system's capability to fully accommodate variable wind generation. After constructing the uncertainty set of wind generation as a polyhedron according to the given probability, the estimation problem is formulated as a three-level optimization problem, which is solved by the presented bisection algorithm. The case study on the IEEE-39 bus system shows that there is relatively low guaranteed chance that power system can fully accommodate variable wind generation under the UC strategy made by the traditional deterministic method. Moreover, it is necessary to introduce novel methods to the decision-making process to cope with wind generation uncertainty. Comparing the capabilities to fully accommodate variable wind generation of strategies made by different methods, studying the impacts of different parameters on the power system's performance may be included in the future work. 


\end{document}